\documentclass[11pt]{article}
\usepackage[ansinew]{inputenc}
\usepackage{graphicx}
\usepackage{amsmath,amssymb}
\usepackage{multirow} 
\usepackage{algorithm,algpseudocode}
\usepackage{subfig}

\newcommand{\Ord}{\mathcal{O}}
\newcommand{\du}{\,\mathrm{d}}
\newcommand{\R}{\mathbb{R}}
\newcommand{\N}{\mathbb{N}}

\newcommand{\e}{\mathrm{e}}

\newcommand{\diag}{\mathrm{diag}}
\newcommand{\abs}{\mathrm{abs}}
\newcommand{\xb}{\mathbf{x}}
\newcommand{\vb}{\mathbf{v}}

\newtheorem{remark}{Remark}

\title{A semi-Lagrangian Vlasov solver in tensor train format} 

\author{Katharina Kormann\thanks{Lehrstuhl f{\"u}r Numerische Methoden in der Plasmaphysik, Technische Universit{\"a}t M{\"u}nchen, Zentrum Mathematik, Boltzmannstr. 3, 85747 Garching, Germany. 
({katharina.kormann@tum.de}). Max--Planck--Institut f{\"u}r Plasmaphysik, Boltzmannstr. 2, 85748 Garching, Germany.}}

\date{}

\begin{document}
\maketitle

\begin{abstract}
In this article, we derive a semi-Lagrangian scheme for the solution of the Vlasov equation represented as a low-parametric tensor. Grid-based methods for the Vlasov equation have been shown to give accurate results but their use has mostly been limited to simulations in two dimensional phase space due to extensive memory requirements in higher dimensions. Compression of the solution via high-order singular value decomposition can help in reducing the storage requirements and the tensor train (TT) format provides efficient basic linear algebra routines for low-rank representations of tensors. In this paper, we develop interpolation formulas for a semi-Lagrangian solver in TT format. 
In order to efficiently implement the method, we propose a compression of the matrix representing the interpolation step and an efficient implementation of the Hadamard product. We show numerical simulations for standard test cases in two, four and six dimensional phase space. Depending on the test case, the memory requirements reduce by a factor $10^2-10^3$ in four and a factor $10^5-10^6$ in six dimensions compared to the full-grid method.

\end{abstract}

\section{Introduction}

The Vlasov equation models the evolution of a plasma in an external or self-consistent field. In its full generality, the model consists of an advection equation in the six-dimensional phase space coupled to Maxwell's equations. Since analytical solutions are usually not known, the numerical simulation of these problems is of fundamental importance. Due to the high dimensionality and the development of small structures the numerical solution is very challenging. There are essentially three classes of solvers that are used in simulations: particle-in-cell (PIC) methods, Eulerian solvers, and semi-Lagrangian methods. The idea of particle methods (cf. e.g. \cite{Birdsall91}) is to distribute a number of macro-particles in the computational domain that are evolved according to the equations of motion. Eulerian as well as semi-Lagrangian solvers, on the other hand, are based on a phase-space grid. In a Eulerian solver (cf. \cite{Arber02} and references therein), the spatial derivatives are approximated by a standard method (e.g. finite volumes \cite{Qiu11} or discontinuous Galerkin \cite{Ayuso12,Cheng13}) and the system is advanced in time using an ODE solver (e.g. Runge--Kutta). Semi-Lagrangian solvers (cf. e.g. \cite{Sonnendruecker99}) update the values of the grid point by evolution along characteristics. They have the advantage that they do not suffer from time step restrictions by the Courant--Friedrichs--Lewy (CFL) condition.

Grid-based methods suffer from the curse of dimensionality, i.e., from the fact that the number of unknows grows with the number of dimensions of the problem. For this reason, grid-based simulations of the 6D Vlasov equations are rarely done. Monte-Carlo methods are computationally less expensive for high-dimensional problems but suffer from a numerical noise problem.

To alleviate the curse of dimensionality for grid-based solvers, several methods especially suited for high-dimensional problems have been developed in the numerical community.  Such methods are the sparse-grid method \cite{Bungartz04}, tensor-based methods \cite{Hackbusch12,Khoromskij12,Grasedyck13}, and radial basis functions \cite{Buhmann00}. In this paper, we consider the solution of the Vlasov--Poisson system in tensor train format \cite{Oseledets11} which is a type of tensor-based methods with attractive numerical properties. The main concept of this method is to represent the solution as nested sums of tensor products. The compression of full-grid data to a sum of tensor products can be achieved by performing high-order singular value decompositions (HOSVD). The HOSVD can also be used to recompress data during time-dependent simulations. Furthermore, all basic numerical operations are defined in the tensor train (TT) format. Even though small filaments evolve in the solution of the Vlasov equation, it has been shown in \cite{Hatch12} that the data from a gyrokinetic Vlasov simulation can be compressed by HOSVD.

Tensor-based methods were introduced in the quantum chemistry community by Meyer et al. \cite{Meyer90} and have been further developed since then. In recent years, mathematical aspects of tensor-based methods have been addressed and formalized in the tensor train \cite{Oseledets12} as well as the hierarchical tensor \cite{Hackbusch09} format. The tensor train format has beed applied to the solution of various high-dimensional differential equations \cite{Kazeev13, Kazeev12, Dolgov12}. In particular, Dolgov et al. \cite{Dolgov14} have designed a tensor train algorithm to solve the Farley--Buneman instability in ionospheric plasma. The equations are similar to the Vlasov--Poisson equations considered in this paper. 

The outline of the paper is as follows. In the next section, we introduce the split-step semi-Lagrangian solver for the Vlasov--Poisson equation. Then we give a short summary of the tensor train format in \S \ref{sec:tt}. The tensor train semi-Lagrangian method is devised in \S \ref{sec:ttsl} and its efficient implementation is discussed in \S \ref{sec:implement}. In particular, we discuss an efficient implementation of the TT Hadamard product and compression of the interpolation operator. Numerical results are provided in \S \ref{sec:numerics} and conclusions as well as future research directions are given in \S \ref{sec:conclusions}.

\section{Vlasov--Poisson equation and semi-Lagrangian discretization} \label{sec:vlasov}

The evolution of the distribution function $f(x,v)$ of a plasma is governed by the Vlasov equation
\begin{equation*}
\partial_t f + v \cdot \nabla_x f + F (t, x, v) \cdot \nabla_v f = 0,
\end{equation*} 
where the force field $F$ is given by the Lorentz force due to external and self-consistent electromagnetic fields. The self-consistent fields can be computed by Maxwell's equations. If the magnetic field is small, it suffices to consider $F(t,x,v) = -E(t,x)$. Then the self-consistent part of the electric field can be computed by the Poisson equation
\begin{equation*}
-\Delta \phi(x,t) = 1-\rho(x,t), \quad E(x,t) = -\nabla \phi(x,t),
\end{equation*}
where $\rho = \int f(x,v) \du v$ is the particle density. In this paper, we focus on such Vlasov--Poisson equations.

The Vlasov equation is a hyperbolic equation and the associated characteristics satisfy the differential system
\begin{equation*}
\frac{\du X}{\du t} = V, \quad \frac{\du V}{\du t} = -F(X,V,t).
\end{equation*}
If the evolution of the characteristic curves is known, the distribution function at time $t$ can be computed from a given initial distribution $f_0$ at time $t=0$ as
\begin{equation*}
f(x,v,t) = f_0(x-X(t),v-V(t)).
\end{equation*}
The idea of the semi-Lagrangian method is to discretize the initial distribution on a mesh. In each time step, the characteristic equations are solved numerically backwards in time and the new solution at the grid points is given by the (interpolated) value of the previous solution at the origin of the corresponding characteristic.

For the Vlasov--Poisson equation, the coefficient of the $x$ gradient is only dependent on $v$ and vice versa. Therefore, a split step method can be designed where (constant) advection along one coordinate direction is considered at a time. In this case, the equations for the characteristics can be solved exactly. This yields the split-step semi-Lagrangian scheme shown in Algorithm \ref{alg:sssl} which was originally introduced by Chen and Knorr \cite{Cheng76}. Note that the interpolation along the $d/2$-dimensional $x$- and $v$-coordinates can be split into $d/2$ separate one-dimensional interpolations.  Various aspects of the semi-Lagrangian method for the Vlasov equations have for instance been discussed in \cite{Filbet01,Besse04,Crouseilles08a,Crouseilles08b,Crouseilles10,Mehrenberger13}.

\begin{algorithm}
\caption{Split-step semi-Lagrangian scheme by Chen and Knorr \cite{Cheng76}}\label{alg:sssl}
Given $f^{(m)}$ and $E^{(m)}$ at time $t_m$, we compute $f^{(m+1)}$ at time $t_m + \Delta t$ as follows:

\begin{enumerate}
  \item Solve $f_t - E^n f_v = 0$ on half time step: $f^{(m,*)}(x_i,v_j) = f^{(m)}(x_i,v_j+E_i^{(m)}\frac{\Delta t}{2})$.
  \item Solve $f_t + v f_x = 0$ on half time step: $f^{(m,**)}(x_i,v_j) = f^{(m,*)}(x_i-v_j\frac{\Delta t}{2},v_j)$.
    \item Compute $\rho(x_i,v_i)$ and solve the Poisson equation for $E^{(m+1)}$.
     \item Solve $f_t - E^{(m+1)} f_v = 0$ on half time step: $f^{(m+1)}(x_i,v_j) = f^{(m,**)}(x_i,v_j+E_i^{(m+1)}\frac{\Delta t}{2})$.
\end{enumerate}
\end{algorithm}

\section{The tensor train format} \label{sec:tt}

On a tensor product grid, the number of points grows exponentially in the dimension variable. In many cases, however, the complete information---or at least an accurate approximation---can be recovered from much less data. The simplest example is a function of the form
\begin{equation*}
f(x_1, \ldots, x_d) = \prod_{k=1}^d f_k(x_k).
\end{equation*}
On a grid with $n_k$ points along dimension $k$, it suffices to store the $\sum_{k=1}^d n_k$ function values $f_k(x_{k,j})$, $j=1, \ldots, n_k$, $k=1,\ldots,d$. The value at any grid point can be reconstructed from this data. Generalizing from this very special case, the tensor train (TT) format \cite{Oseledets11} offers the possibility of representing a multidimensional function as nested sums of such Kronecker products which yield good approximations of much more complicated functions. In the tensor train format, each dimension $k$ is represented by a kernel $Q_{k} \in \R^{r_{k-1} \times n_k \times r_k}$, in which the second index runs over the grid points along the $k$th dimension and the first und third index take care of couplings to the $(k-1)$th and $(k+1)$th dimension, respectively. The value at the grid point $(x_{i_1},\ldots,x_{i_d})$ can be reconstructed as
\begin{equation*}
f(x_{i_1},\ldots,x_{i_d}) = \sum_{\alpha_0=1}^{r_0}\ldots\sum_{\alpha_d=1}^{r_d} \prod_{k=1}^d Q_k(\alpha_{k-1},i_k,\alpha_k).
\end{equation*}
The size of the \emph{ranks} $r_k$, $k=1,\ldots,d-1$, depends on the structure of the function and the level of accuracy required. Since $f$ is a scalar function, we have $r_0 = r_d = 1$ and we will sometimes omit the corresponding index.

In case we have a function represented as a $d$-dimensional tensor, we can use a sequence of singular value decompositions (SVD) applied to matrifications of the tensor to find a representation in TT format to a given error tolerance or with a fixed maximum rank.

The tensor train format offers not only compression of high-dimensional data but also provides simple algorithms for basic tensor-tensor and matrix-tensor operations. For instance, we can build the kernels of the sum of two TT tensors $A=Q_1^{A} \cdot \ldots \cdot Q_d^A$ and $B=Q_1^{B} \cdot \ldots \cdot Q_d^B$ by setting
\begin{equation*}\begin{aligned}
Q_1^{A+B}(i_1) &= \begin{pmatrix}
Q_1^{A}(i_1) & Q_1^B(i_1)
\end{pmatrix}, \quad Q_d^{A+B}(i_d) = \begin{pmatrix}
Q_d^{A}(i_d) \\ Q_d^B(i_d)
\end{pmatrix}, \\
 Q_k^{A+B}(i_k) &= \begin{pmatrix}
Q_k^{A}(i_k) & 0 \\ 0 & Q_k^B(i_k)
\end{pmatrix}, k=2,\ldots,d-1.
\end{aligned}\end{equation*}
This operation is essentially a copying function. However, the ranks of $A$ and $B$ sum up to the ranks of $A+B$. In a matrix-vector product, the ranks of the matrix and the vector are even multiplied. Obviously, continued application of basic operations, for instance to propagate the tensor in time, will destroy the compression. Most often, however, the representation of the new TT tensor can be truncated and it is essential to continuously add rounding steps to any algorithm operating on TT tensors. Since one wants to be able to truncate one rank at the time, a left-to-right sweep with QR-decompositions of the kernels is performed to orthonormalize all kernels except for the last. Then, the kernels are singular-value decomposed individually in a right-to-left sweep where only the non-orthogonal kernel is touched in each iteration. Algorithm \ref{alg:ttround} implements the TT rounding. Note that we use an \emph{absolute} threshold in constrast to the rounding in \cite{Oseledets11}. The rounding requires the computations of $d-1$ QR decompositions (for orthonormalization) as well as $d-1$ SVD. Since we are not interested in the zero blocks, it suffices to compute an economy-size decomposition in both cases, i.e., the vectors corresponding to zero blocks are left out. Computing such economy-size QR or SV decompositions for a $r_{k-1}n \times r_k$ matrix is of complexity $\Ord(r_{k-1}nr_k^2)$ (cf. \cite[Chapt. 2.5]{Hackbusch12}).
 The complexity of a rounding step is hence $\Ord(dnr^3)$ where we have assumed all ranks to be equal to $r$ and $n$ grid points along each dimension. A more detailed description of operations in TT format can be found in \cite{Oseledets11}.

\begin{algorithm}[H]
\caption{Tensor train rounding (from \cite{Oseledets11}).}\label{alg:ttround} 
\begin{algorithmic}
\Require{$d$-dimensional tensor in TT format $A=Q_1 \cdots Q_d$; absolute tolerance $\varepsilon$ and maximum rank $r_{\max}$}
 \Ensure{$d$-dimensional tensor $B$ in TT format with kernels rounded according to input tolerance $\varepsilon$ and/or ranks bounded by $r_{\max}$}
 \State{ $\delta = \frac{\varepsilon}{\sqrt{d-1}}$}\Comment{Scale truncation parameter} 
 \For{k=1 to d--1}
 	\State{[$Q_k(\beta_{k-1}i_k ,\beta_k),R(\beta_{k},\alpha_k)$] = QR($Q_k(\beta_{k-1}i_k,\alpha_k)$)}
 	\State{$Q_{k+1}(\beta_{k},i_{k+1}\alpha_{k+1}) = R(\beta_k,\alpha_k) Q_{k+1}(\alpha_{k},i_{k+1}\alpha_{k+1})$}
 \EndFor
 \For{k=d to 2}
 \State{[$U(\beta_{k-1},\gamma_{k-1})$,$S$,$Q_k(\gamma_{k-1},i_k\gamma_k)^T$] =  SVD$_{\delta,r_{\max}}$($Q_k(\beta_{k-1},i_k\gamma_{k})$)} \Comment{$\delta$- truncated SVD with maximum rank $r_{max}$}
 \State{$Q_{k-1}(\beta_{k-2}i_k,\gamma_{k}) =  Q_{k-1}(\beta_{k-2}i_k,\beta_{k-1}) U(\beta_{k-1},\gamma_{k-1}) S$ }
 \EndFor
 \end{algorithmic}
\end{algorithm}

\section{A semi-Lagrangian method in tensor train format}\label{sec:ttsl}

In this section, we explain how a split-step semi-Lagrangian method can be designed in tensor train format. First, we derive the formulas for the example of linear interpolation in two dimensions (i.e. 1D Vlasov) before discussing other interpolation formulas and higher dimensions. We also discuss suitable ordering of the coordinates in four and six dimensions and the effects of TT rounding. The conservation properties of the method are discussed and we propose a projection to the manifold spanned by constant mass and momentum. Finally, we shortly discuss the solution of the Poisson problem.

\subsection{Derivation for 1D Vlasov}\label{sec:ttsl2D}

Consider the $x$-advection in two dimensions. We assume a tensor train representation of the distribution function at time $t_m$ of the form
\begin{equation*}
 f^{(m)}(x,v) \approx \sum_{\alpha} Q_1^{(m)}(x,\alpha) Q_2^{(m)}(\alpha,v).
\end{equation*}
We now consider the displacement in $x$ direction by $-\Delta t v$. To keep the presentation simple, we derive the formulas for linear interpolation. Even though not necessary for stability, we impose the CFL-like condition that
\begin{equation}\label{eq:sl_cfl}
  |\Delta t v| \leq \Delta x,
\end{equation}
where $\Delta x$ is the grid spacing along $x$.
On a full grid, the distribution function at the new time step would be computed according to the formula
\begin{equation*}\begin{aligned}
  f^{(m+1)}(x_j,v_k) &= f^{(m)}(x_j-\Delta t v_k,v_k) \approx \max\left(0,\frac{\Delta t v_k}{\Delta x}\right) f^{(m)}(x_{j-1},v_k) \\&+ \left(1-\abs\left(-\frac{\Delta t v_k}{\Delta x}\right)\right)f^{(m)}(x_j,v_k) + \max\left(0,-\frac{\Delta t v_k}{\Delta x}\right) f^{(m)}(x_{j+1},v_k),
\end{aligned}\end{equation*}
where the indices are periodically shifted for periodic boundary conditions. Note that the displacement can either be positive or negative. For a function in tensor train format the kernels representing $f^{(m+1)}$ are computed from the kernels of $f^{(m)}$ by

\begin{equation}\label{eq:sl_shift}\begin{aligned}
  &Q_1^{(m+1)}(x_j,\alpha) Q_2^{(m+1)}(\alpha,v_k) = Q_1^{(m)}(x_j-\Delta t v_k,\alpha) Q_2^{(m)}(\alpha,v_k) \\
\quad &\approx Q_1^{(m)}(x_{j-1},\alpha) \left( Q_2^{(m)}(\alpha,v_k) \max\left(0,\frac{\Delta t v_k}{\Delta x}\right)\right)\\
&\quad +  Q_1^{(m)}(x_{j},\alpha) \left(Q_2^{(m)}(\alpha,v_k)\left(1-\abs\left(-\frac{\Delta t v_k}{\Delta x}\right)\right)\right) \\
&\quad +  Q_1^{(m)}(x_{j+1},\alpha) \left(Q_2^{(m)}(\alpha,v_k)\max\left(0,-\frac{\Delta t v_k}{\Delta x}\right)\right).
\end{aligned}\end{equation}
This can be interpreted as the sum of three tensor trains. In each case, the first kernel is a shifted version of the original kernel und the second kernel is scaled depending on the value of $v$. 
Hence, we can compute the advection in $x$ direction performing the following steps:
\begin{enumerate}
  \item Form the three shifted kernels of $Q_1^{(m)}$.
  \item Form three scaled $Q_2^{(m)}$-kernels that are line-wise multiplied by $\max\left(0,\frac{\Delta t v_k}{\Delta x}\right)$, $\left(1-\abs\left(-\frac{\Delta t v_k}{\Delta x}\right)\right)$, and $\max\left(0,-\frac{\Delta t v_k}{\Delta x}\right)$, respectively.
  \item Form three TT-tensors from the resulting kernels.
  \item Add the TT-tensors and perform a rounding step.
\end{enumerate}

Each of the tensors formed in step 1 and 2 have the same rank as the original tensor. Adding the tensors will increase the rank (by a factor three in this case). However, the rank can usually be reduced again by performing a rounding step. 

We can also write the evolution operator as a matrix $A$ in tensor product form. If we denote by $S_{j}$ the matrix with one diagonal of ones shifted by $j$ from the center, we have
\begin{equation*}\begin{aligned}
A =& S_{-1} \otimes \diag\left(\max\left(0,\frac{\Delta t v_k}{\Delta x}\right)\right) + S_0 \otimes \diag\left(1-\abs\left(-\frac{\Delta t v_k}{\Delta x}\right)\right) \\
&+ S_1 \otimes \diag\left(\max\left(0,-\frac{\Delta t v_k}{\Delta x}\right)\right).
\end{aligned}\end{equation*}
This can be written in TT format as a matrix with rank $r_1=3$.
It is obvious that the advection with respect to $v$ can be done in the same way by interchanging the roles of $Q_1$ and $Q_2$.

In principle, we can use any other interpolator in our derivation. Especially, a centered Lagrange interpolator that includes $p$ points will result in $p+1$ TT-tensors that have to be formed by shifting the $Q_1$-kernel, scaling the $Q_2$-kernel, and finally adding the $(p+1)$ kernels. For a non-nodal interpolator, like splines, step 2 needs to be augmented. Before shifting the kernel, the values of the interpolator weights have to be computed for each column $Q_1^{(n)}(:,\alpha)$.

\begin{remark}
Similar to a Eulerian solver and opposed to the usual case for semi-Lagrangian solvers we have imposed the CFL-like condition \eqref{eq:sl_cfl}. However, it is possible to relax this condition. The condition was not introduced to ensure stability but to make sure that we only have to consider $p+1$ points for a centered interpolator with $p$ points. We can relax the condition at the price of additional terms in the sum \eqref{eq:sl_cfl}. For the condition
\begin{equation}\label{eq:sl_cfl2}
 |\Delta t v| \leq m \Delta x
\end{equation}
with some $m \in \N$, the number of points will be $p+2m-1$. Since $p$ intervals will be used at a time, the scaled $Q_2$-tensors will contain an increasing number of zeros (as $m$ increases). This might be exploited to further improve on the efficiency. 
\end{remark}

\subsection{Generalization to higher dimensions}

In higher dimensions, the 2D algorithm can be applied in essentially the same way to parts of the kernels while others are left unchanged. In particular, the advection along a spatial dimension will only depend on one (velocity) dimension also in 4D and 6D. Then, we treat the two corresponding kernels as discussed in the previous section and keep the other two or four kernels unchanged.

For the velocity advection, the situation becomes more involved. For simplicity, we consider the 4D case. The generalization to 6D is straight forward. The displacement is now not only dependent on one dimension but on two.  
Let us revisit the linear interpolation and consider the advection along $v_1$
\begin{equation}\label{eq:sl4d_shift}\begin{aligned}
  &Q_1^{(m+1)}(x_{1,j_1},\alpha_1) Q_2^{(m+1)}(\alpha_1,x_{2,j_2},\alpha_2)Q_3^{(m+1)}(\alpha_2,v_{1,j_3},\alpha_3) Q_4^{(m+1)}(\alpha_3,v_{2,j_4}) =  \\
& Q_1^{(m)}(x_{1,j_1},\alpha_1) Q_2^{(m)}(\alpha_1,x_{2,j_2},\alpha_2)Q_3^{(m)}(\alpha_2,v_{1,j_3}+\Delta t E_1(x_{1,j_1},x_{2,j_2}),\alpha_3) Q_4^{(m)}(\alpha_3,v_{2,j_4}).
\end{aligned}\end{equation}
The displacement $\Delta t E_1(x_{1,j_1},x_{2,j_2})$ is a function of two variables and we assume we are given its TT representation. For the linear interpolation, we need a TT representation of $g_1(x_1,x_2) = \max\left(0,-\frac{\Delta t E_1(x_1,x_2)}{\Delta v}\right)$, $g_2(x_1,x_2)=1-\abs\left(\frac{\Delta t E_1(x_1,x_2)}{\Delta v}\right)$, and $g_3(x_1,x_2) = \max\left(0,\frac{\Delta t E_1(x_1,x_2)}{\Delta v}\right)$. Let us denote them by $g_i(x_1,x_2) = W^{(i)}_1(x_1,\gamma^{(i)})  W^{(i)}_2(\gamma^{(i)},x_2)$ and by $s^i$ the ranks of each TT tensor. Then \eqref{eq:sl4d_shift} becomes
\begin{equation*}
\begin{aligned}
  &\sum_{i=1}^3\left( \sum_{\gamma^{i}=1}^{s^{i}}\left( Q_1^{(m)}(x_{1,j_1},\alpha_1) W_1^i(x_{1,j_1},\gamma^i) \right) \left(Q_2^{(m)}(\alpha_1,x_{2,j_2},\alpha_2)W_2^i(\gamma^i,x_{2,j_2})\right)\right) \cdot \\
&\quad Q_3^{(m)}(\alpha_2,v_{1,j_3-2+i},\alpha_3) Q_4^{(m)}(\alpha_3,v_{2,j_4}).
\end{aligned}\end{equation*}
Thus, we have to form (for each $i$) $s^i$ scaled versions of the 2D TT tensor represented by $Q_1,Q_2$ and add them up. To keep the size of the rank small, it is advantageous to truncate after each addition. Note that this operation can also be described as a multiplication of the TT-tensor $Q_1 Q_2$ by the matrix $\diag\left(W_1^{(i)}W_2^{(i)}\right)$. This is, of course, much more expensive than step 3 of the 2D algorithm but there is no additional difficulty due to the use of the TT format. 

Now, we turn to the question of how to compute the kernels $W^{(i)}_1W^{(i)}_2$. 
For this, we need to compute the positive and negative part of a TT tensor. This cannot easily be done because the value at the various grid points is not stored explicitly in TT format. 
 Since the displacement is just a 2D (or 3D) object, one might as well compute its values on the full grid to perform these operations. However, the positive and negative part are non-smooth functions and the compression in TT format will generally be rather poor. 
 
This problem is not specific to linear interpolation. The important fact is that we use a different interpolation function depending on the interval into which the point is displaced.  An alternative is to always use the Lagrange polynomial computed on an odd number of points around the original point.  As long as we impose a CFL-like condition that makes sure that we do not displace more than the grid size, we interpolate close to the center. In this case, we only have to compute polynomials of the displacement which can easily be done in TT format. However, we have to make sure the displacements stay small, i.e., relaxing the CFL-like condition will generally result in rather poor approximations. However, the displacement in the $v$-advection step due to the electric field is usually small compared to the displacement in the $x$-advection step. Therefore, relaxing the CFL-like condition for the $x$-advection step only will already result in a fairly efficient time stepping.

\subsection{Ordering of variables}

The natural ordering of the coordinates is to start with the spatial coordinates and then add the velocity coordinates. On a full grid, a reordering does not change the algorithm. In the TT representation, however, we have an explicit coupling between neighboring dimensions. Hence, the compression is affected by the ordering of the variables. To illustrate this, we consider the three-variate function $f(x,y,z) = (f_1(x)+f_3(z))f_2(y)$. To represent this function, we need a TT tensor with ranks $r_1=r_2 = 2$ and kernels
\begin{equation*}\begin{aligned}
&Q_1(i_1,1) = f_1(x_{i_1}), \quad Q_1(i_1,2) = 1,\\	
&	Q_2(1,i_2,1) = f_2(y_{i_2}), \quad Q_2(1,i_2,2) = 0, \quad Q_2(2,i_2,1) = 0, \quad Q_2(2,i_2,2) = f_2(y_{i_2}), \\
&	Q_3(1,i_3) =1, \quad Q_3(2,i_3) =  f_3(z_{i_3}).
\end{aligned}	\end{equation*}
If we instead reorder the variables as $x,z,y$, the TT tensor representing $f$ has ranks $r_1=2$ but $r_2=1$ and kernels
\begin{equation*}\begin{aligned}
&Q_1(i_1,1) = f_1(x_{i_1}), \quad Q_1(i_1,2) = 1,\\	
	&Q_2(1,i_2,1) =1, \quad Q_2(2,i_2,1) =  f_3(z_{i_2}),\quad
	Q_3(1,i_3) = f_2(y_{i_3}).\\
\end{aligned}\end{equation*}

Analyzing the splitting algorithm, we see that the $x_i$-advection step couples dimensions $x_i$ and $v_i$ and the $v_i$ advection step couples $v_i$ and $x_1,\ldots, x_{d/2}$. Therefore, it is reasonable to assume that an ordering that groups the pairs $(x_i,v_i)$ as well as the spatial variables together will result in better compression. In four dimensions, a reordering of the coordinates as $v_1$, $x_1$, $x_2$, $v_2$ satisfies all the requirements. Also, it simplifies the advection steps which---up to an initial orthogonalization steps---only involves two or three neighboring kernels of the TT tensor. Moreover, we note that $r_0=r_4=1$ which is why the first and last kernels are only 2-tensors and therefore generally smaller than the inner kernels that are 3-tensors. This further improves the compression if the number of grid points along the velocity dimensions is higher than along the spatial dimensions. This situation is not uncommon in Vlasov--Poisson simulations.

For the six dimensional case, there is no ordering that places together all different coordinate combinations for the six advection steps. Since the coupling appears to be strongest between the pairs $(x_i,v_i)$, we have found the ordering $v_1$, $x_1$, $x_2$, $v_2$, $x_3$, $v_3$ to be most efficient (up to index shifting).
 
\subsection{TT rounding}

When simulating the Vlasov equation with the semi-La\-grangian split-step method in TT format, we constantly compress the data for the distribution function. This adds to the numerical error of the method. In each simulation, we have to decide when to truncate the HOSVD computed to recompress the data. One strategy would be to keep all the information up to round-off errors. Alternatively, we could choose the drop tolerance such that the error from TT rounding is on the same order of magnitude as other numerical errors. Even though the first strategy has its advantages, especially with respect to the conservation properties as discussed in the next section, this will generally become rather expensive. The reason is that the tensors have to resolve numerical errors that are generally much less smooth than the actual solution. Hence, we will need considerably larger ranks in order to resolve spurious information. It is therefore recommended to choose the drop tolerance carefully to fit the numerical errors of the underlying method. However, we need an error estimator for the underlying spatial discretization in order to be able to automatically decide on a proper tolerance.

When simulating the Vlasov equation over longer times, filaments evolve. This means that the distribution function is relatively well-resolved in the beginning. A simple strategy to account for this fact is to scale the tolerance $\varepsilon$ at final time according to the time step $j$ as  $\frac{j}{N_t}\varepsilon$, where $N_t$ is the total number of time steps in the simulation. We have used this scaling in our numerical experiments.

\subsection{Conservation properties} \label{sec:sltt_conservation}

Many integrals of the solution of the Vlasov equation are conserved: mass, momentum, energy, and all $L_p$ norms. Moreover, the maximum and minimum value of the solution remain constant. These invariants of the mathematical function mimic the physical behaviour of the distribution function. Therefore, a good numerical method should conserve these properties as accurately as possible. Without rounding in the TT representation, we would---up to roundoff errors---recover the same solution as the underlying method on the full grid. Hence, our method would inherit the conservation properties of the underlying full-grid method. If we perform rounding up to some drop tolerance, we fulfill the conservation laws of the full-grid method with the accuracy of the rounding. The conservation thus depends on the drop tolerance. When choosing a loose tolerance, a
 projection onto the manifold spanned by one or more conservation laws should be considered (cf. the next section). Since the SVD yields a best approximation in $\ell_2$ sense, the $\ell_2$ norm is expected to be conserved to a large extent.

\subsection{Projection}

If the tolerance is chosen too loose or fixed ranks are used, one can use a projection method (cf. \cite[Chapt. IV.4]{Hairer06}) to improve conservation. Given the propagated solution, the closest solution on the manifold defined by the conserved quantities is found by minimizing the Lagrange function describing this constrained minimization problem.  In particular, we consider conservation of mass and momentum. The discrete mass and momentum are defined as
\begin{equation*}\begin{aligned}
&\text{mass}(f) = (\Delta x)^d(\Delta v)^d \sum_{(\xb,\vb) \in \mathcal{G}} f(\xb,\vb)\\
&\text{mom}(f) = (\Delta x)^d(\Delta v)^d \sum_{(\xb,\vb) \in \mathcal{G}} f(\xb,\vb) \vb,
\end{aligned}\end{equation*}
where $\mathcal{G}$ denotes the set of grid points. The number of points is denoted by $|\mathcal{G}|$. Given the solution $\tilde f^{(m)}$ at point $(\xb,\vb)$, obtained by the time evolution algorithm, we add a perturbation by $(1, \vb^T)\mu$ for a suitably chosen Lagrange multiplier $\mu \in \R^{d+1}$. In the case of mass and momentum conservation, the projection is rather simple since all the $d+1$ projections are orthogonal to each other. The projected solution is then given by
\begin{equation*}
f^{(m)}(\xb,\vb) = \tilde f^{(m)} + \frac{\text{mass}(f^{(0)})-\text{mass}(\tilde f^{(m)})}{(\Delta x)^d(\Delta v)^d|\mathcal{G}|} + \sum_{k=1}^d\frac{\text{mom}_k(f^{(0)})-\text{mom}_k(\tilde f^{(m)})}{(\Delta x)^d(\Delta v)^d|\mathcal{G}|}  v_k.
\end{equation*}
Applying this projection to the solution in TT format will increase the rank by $d+1$ since we add $d+1$ rank-one tensors to the solution.

\subsection{Poisson's equation}

So far, we have only discussed the solution of the Vlasov equation. However, we also have to solve a Poisson problem in each time step. Since the Poisson equation only depends on the spatial variable, the dimensionality is cut into half. In our prototype implementation, we have therefore chosen to solve Poisson's equation with a pseudo-spectral method on the full grid. However, a pseudo-spectral solver based on the fast Fourier transform in TT format \cite{Dolgov12a} will presumably improve the efficiency of our method.

\section{Efficient implementation}\label{sec:implement}

While it is possible to achieve considerable data compression when using the TT format, it is less obvious if this helps in reducing the computing time. Even though a much smaller amount of data needs to be handled, we have to rely on more complex algorithms.
As mentioned in \S \ref{sec:tt} the ranks are multiplied in a matrix-vector product in TT format and we have to perform a rounding step together with each matrix-vector product. In this section, we discuss the TT matrix-vector product in more detail and explain how to efficiently implement it for the matrices appearing in our semi-Lagrangian solver.

\subsection{The TT matrix-vector product}

Let us consider the matrix-vector product $w=Mu$, where all objects are in TT format. When counting arithmetic operations, we assume all ranks of the matrix to be $s$ and all ranks of the vector to be $r$ and the matrix to be quadratic. In practice, the complexity will be dominated by the maximum rank. 

The complexity of  a direct matrix vector product, i.e., multiplying the kernel and then performing a TT rounding on $w$ at the end, is $\Ord(d n^2 r^2 s^2 + d n r^3 s^3)$, where the first part is due to the matrix vector product and the second due to the TT rounding. If the matrix is a (nested sum of) Kronecker products of sparse matrices in each dimension, the complexity is reduced to $\Ord(d n r^2 s^2 + d n r^3 s^3)$. In this case, the complexity is clearly dominated by the rounding operation.  Already if the ranks are on the order 10, the factor $r^3s^3$ becomes significant. Compared to the complexity of the sparse matrix-vector product on the full grid, $\Ord(n^d)$, this might not be small for $d=2,4,6$. 
However, the optimal rank of $w$ is usually close to $\max(r,s)$ rather than $rs$. Therefore, the complexity can be improved if the matrix-vector product is not computed directly but approximation to a certain threshold or maximum rank size is incorporated into the matrix vector product. 

As an alternative to SVD-based rounding, methods within the alternating direction framework pose the problem of finding a low-rank approximation to a tensor as an optimization problem.  These methods are iterative and optimize on one (or two) kernels at a time while the others are frozen. The ALS method \cite{Lathauwer00} works with fixed ranks and has a complexity of $\Ord(dnsr^3+dn^2s^2r^2+dnr^3)$ for full or $\Ord(dnsr^3+dns^2r^2+dnr^3)$ for sparse matrix kernels.  However, the convergence is rather slow and there are no convergence estimates. The DMRG method or MALS \cite{Holtz12,White92} is an alternative that is based on the same algorithmic idea but operates on two kernels at a time. In this way, the ranks can be adapted and convergence is faster and better understood. However, the complexity is increased by a factor $n$. Matrix-vector products based on the DMRG method were introduced in \cite{Oseledets11a}.  Recently, the AMEn routine \cite{Dolgov13,Dolgov13a} has been presented that allows for adaptive ranks and relatively fast convergence with complexity comparable to the ALS method.

\subsection{Efficient multiplication with a diagonal matrix -- Hadamard product}\label{sec:diagonal_matrix}

A standard TT representation of a matrix does not take sparsity of the one-dimensional kernels into account. The multiplication by a diagonal matrix should therefore be considered as a Hadamard product, i.e., the element-wise product of two tensors.

As seen in Algorithm \ref{alg:ttround}, the TT-rounding algorithm proceeds in two steps:
\begin{enumerate}
	\item Left-to-right sweep with orthogonalization of the kernels $1, \ldots, d-1$ via QR decomposition.
	\item Right-to-left sweep with SVD and truncation of singular values.
\end{enumerate}

For the first step, we note that we can find an orthogonalized representation of both TT tensors. If we then compute the Hadamard product of the kernels, the resulting TT tensor is again orthogonal. This splitting of the orthogonalization step reduces the complexity from $(d-1)ns^3r^3$ to $(d-1)n(s^3 + r^3)$ compared to when computing the QR decomposition of the multiplied kernel.

For the SVD in the second step, let $\sigma=\min(s,r)$ and $\tau = \max(s,r)$ and consider the kernel in dimension $j$. Since we have already truncated over the $j$th rank, we can assume that this rank is of the order $\tau$ and rank $j-1$ should be truncated. The multiplied kernel consists of $\sigma$ blocks size $\tau$. We now take two such blocks and truncate them. Then, we add another block of size $\tau$ and truncate again. This is repeated until all $\sigma$ blocks are included. The procedure is summarized in Algorithm~\ref{alg:tthadamard}.

In total, the complexity is reduced to the order $\Ord(dn r s \max(r,s)^2)$, if we assume that the rank of the product is approximately $\max(r,s)$. The reduction is  about a factor $\min(r,s)^ 2/8$ compared to the direct method. Of course, one could consider any other grouping of the constituents of the kernel. For instance, one might group more than two $\tau$-sized kernels together.

\begin{algorithm}
 \caption{Rounded tensor train Hadamard product.}\label{alg:tthadamard}
 \begin{algorithmic}
 \Require{$d$-dimensional tensors in TT format $A=Q_1^A \cdots Q_d^A$ and $B=Q_1^B \cdots Q_d^B$; tolerance $\varepsilon$ and maximum rank $r_{\max}$}
 \Ensure{$d$-dimensional tensor $C=Q_1 \cdots Q_d$ in TT format being the Hadamard product $A \star B$ with kernels rounded according to input tolerance and/or ranks bounded by $r_{\max}$}
 \For{k=1 to d--1}\Comment{Orthogonalize kernels of $A$}
 \State{	[$Q_k^A(\beta_{k-1}^Ai_k ,\beta_k^A),R(\beta_{k}^A,\alpha_k^A)$] = QR($Q_k^A(\beta_{k-1}^Ai_k,\alpha_k^A)$)}
 	\State{$Q_{k+1}^A(\beta_{k}^A,i_{k+1}\alpha_{k+1}^A) = R(\beta_k^A,\alpha_k^A) Q_{k+1}(\alpha_{k}^A,i_{k+1}\alpha_{k+1}^A)$}
 \EndFor
  \For{k=1 to d--1}\Comment{Orthogonalize kernels of $B$}
 	\State{[$Q_k^B(\beta_{k-1}^Bi_k ,\beta_k^B),R(\beta_{k}^B,\alpha_k^B)$] = QR($Q_k(^B\beta_{k-1}^Bi_k,\alpha_k^B)$)}
 	\State{$Q_{k+1}^B(\beta_{k}^B,i_{k+1}\alpha_{k+1}^B) = R(\beta_k^B,\alpha_k^B) Q_{k+1}^B(\alpha_{k}^B,i_{k+1}\alpha_{k+1}^B)$}
 \EndFor
 \For{k=1 to d}\Comment{Compute Hadamard product}
 \State{$Q_j = Q_j^A \star Q_j^B$}
 \EndFor
  \State{$\delta = \frac{\varepsilon}{\sqrt{d-1}}$}
  \Comment{Scale truncation parameter}
 \For{k=d to 2}\Comment{Truncate kernels.}
 \State{$\sigma =\min(r_{j-1}^A,r_{j-1}^B)$;  $\tau =\max(r_{j-1}^A,r_{j-1}^B)$}
 \State{$\rho = \min(2,\sigma)\tau$; $iter = \max(\sigma-1,1)$}
 \State{$\tilde \delta = \delta/iter$}
 \For{j=1 to iter}
 \State{[$U(\beta_{k-1},\gamma_{k-1})$,$S$,$Q_k(1:\rho',i_k\gamma_k)^T$] =  SVD$_{\tilde\delta,r_{\max}}$($Q_k(1:\rho,(i_k:)\gamma_{k})$)} \Comment{$\delta$- truncated SVD with maximum rank $r_{max}$}
 \State{$Q_{k-1}(\beta_{k-2}i_k,1:\rho')$ =  $Q_{k-1}(\beta_{k-2}i_k,1:\rho) U(\beta_{k-1},\gamma_{k-1}) S$}
 \State{$\rho = \rho'+\tau$}
 \EndFor
 \EndFor
 \end{algorithmic}
\end{algorithm}

\subsection{Combination of diagonal kernels with one non-diagonal kernel}\label{sec:eff_mvp2}

The idea presented in the previous subsection is limited to diagonal matrices, since the QR decomposition of a sparse matrix is generally non-sparse. Hence, the orthogonal representation of the kernels of the matrix in TT format would be non-sparse, which would give complexities in the range of $n^2$. However, we can essentially apply Algorithm~\ref{alg:tthadamard} if we have a matrix that is non-diagonal in dimension $d$ only since the kernels $Q_d^{A/B}$ are not QR decomposed. We only need to replace the Hadamard product along dimension $d$, $ Q_d^A \star Q_d^B$, by a sparse matrix-vector product.

Using RQ instead of QR decompositions, we can interchange the direction of the loops and create a TT tensor that is non-orthogonal in $Q_1$. Combining QR and RQ decompositions, we can create an orthogonalized TT tensor with any non-orthogonal kernel. In this case, however, we have to choose the direction of the truncation step, i.e., the truncation step will only include the initially non-orthogonalized kernel together with either the kernels with smaller or larger indices. If we do not want to truncate on one side of the non-diagonal kernel---for instance because the matrix is the identity on one side---Algorithm \ref{alg:tthadamard} is still applicable. As long as we order the dimensions such that the spatial coordinates are adjacent, the advection matrices in our splitting semi-Lagrangian scheme have this structure.

\subsection{Combinations of stencils and coefficients through splitting}\label{sec:eff_mvp3}

A more flexible alternative to combine the efficient Hadamard product with a sparse non-dia\-go\-nal kernel is to split diagonal and off-diagonal parts. For the Kronecker product of two matrices $A$ and $B$, it holds that
$$A \otimes B = (I_1 \otimes B) (A \otimes I_2),$$
where $I_{1/2}$ denote identity matrices of the corresponding size. Hence, we can apply all non-diagonal kernels first followed by an application of the diagonal kernels according to Algorithm \ref{alg:tthadamard}. Of course, this means we are applying the matrix in two steps with an error in each step if we apply intermediate rounding in contrast to the alternative discussed in the previous section. On the other hand, the structure of the matrix is not limited. Also, we do not have to explicitly form the kernels of the non-diagonal matrices since no orthogonalization  is necessary. 

\begin{remark}
The leading-order complexity of the various variants of the matrix-vector product that were presented in this section is of the same order as the AMEn matrix-vector product. However, the computation does not involve any iterative method and the constant is therefore supposed to be smaller in general. Indeed, we have seen in numerical experiments that our matrix-vector product is generally faster than the AMEn product. Possibly, the computing time can be further reduced by applying the AMEn algorithm not to the full matrix-vector product but to the rounding steps in the Hadamard-based algorithms.
\end{remark}

\subsection{Rounding of propagation matrices}

When we are explicitly forming the propagation matrices, we may round the matrix before computing the matrix vector product. Especially when using higher order and for the velocity advections in four and six dimension where ranks of the propagation matrices become larger, rounding can reduce the complexity of the advection step.

If we build the full kernels of a matrix, the rounding has a computational complexity of $\Ord(dn^2 s^3)$ for ranks $s$ and $n$ points per dimension. However, there will be many zero entries in this matrix. For the rounding operation, an $m \times n$ TT matrix kernel is treated as a TT tensor kernel of size $m n$. If we have an index $i \in \{1,\ldots,mn\}$ such that the corresponding entries of the kernel are all zero, i.e.,  $Q(:,i,:)=0$, this dimension will always give a zero contribution. Hence, we can erase these dimensions from the TT tensor representing the TT matrix. For a diagonal kernel this means that only the diagonals need to be stored and the size can be reduced from $n^2$ to $n$. For a sparse matrix, we need to keep a sparsity pattern that includes the sparsity pattern of all the sparse matrices   representing the kernel. This reduces the complexity to the order $\Ord(dn s^3)$.

The sparse matrices appearing in our propagation matrices represent an index shift. If we use an interpolator that involves $p$ points around each point, we have $p$ index-shifting matrices and the total sparsity pattern includes a band of $p$ points around the diagonal. However this band of $p$ points will be exactly the same for each of the $n$ points. Since the rank coupling to the neighboring dimensions is the same for all points, it suffices to keep one copy of this band. This reduces the size of the kernel representing the sparse index-shifting matrix from $n^2$ for a full matrix representation to $p$. Hence, the size of the kernel representing the index-shift matrix is independent of the number of grid points.

\subsection{Efficient advection}

In our experiments, we found a matrix-vector product following \S \ref{sec:eff_mvp2} to be most efficient if the advection coefficient only depends on one variable, i.e., velocity advection in two dimensions and all spatial advections. Essentially, this is the algorithm described in \S \ref{sec:ttsl2D} with some specified interpolation formula and applied to parts of an orthonomal TT-tensor. For the case where the coefficient is multivariate, we use the splitting described in \S \ref{sec:eff_mvp3}. In this case, we also round the propagation matrix before computing the matrix-vector product. However, we have observed considerable loss in accuracy when using the same threshold as for the rounding of the TT tensor. Therefore, the threshold is reduced by a factor 4 for the matrix rounding. When rounding the propagation matrix, most redundancies in the matrix-vector product are already eliminated. Hence, the splitting of the SVD over $\sigma$ in Algorithm \ref{alg:tthadamard} does generally not speed up the product.

\section{Numerical results}\label{sec:numerics}

In this section, we present results obtained with the semi-Lagrangian method in TT format for the weak and strong Landau damping as well as the two stream instability. In the simulations, we use a cubic spline interpolator for all advections with univariate coefficient and a five-point Lagrange interpolation for all advections with multivariate coefficients. As a reference, we compare our result to a full grid solution using cubic spline interpolation. All experiments are performed in MATLAB with a prototype implementation based on the TT-Toolbox\footnote{Available at https://github.com/oseledets/TT-Toolbox, downloaded on March 19, 2014.}. The reported computing times are for an Intel Ivy Bridge notebook processor with two cores at 3.0 GHz. 

\subsection{Weak Landau damping}

\begin{figure}[h]
\begin{center}
\includegraphics[scale=0.7]{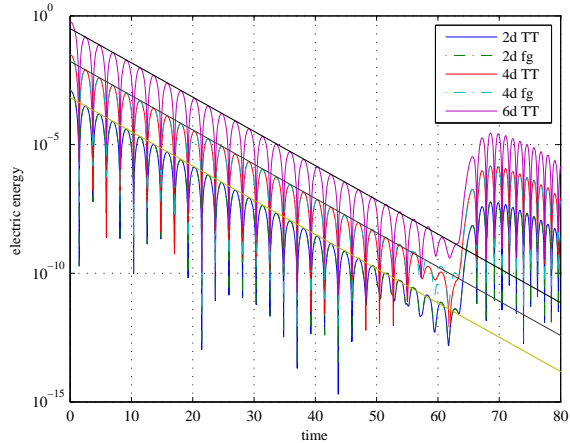}
\caption{Weak Landau damping. Electric energy for simulations on $32^d \times 128^d$ grid.}
\label{fig:linlan_eenergy}
\end{center}
\end{figure}

The initial condition for the standard Landau test case \cite{Krall73} in $d=1,2,3$ dimensions is given by
\begin{equation}\label{eq:weaklan_aligned}
f_0(\xb,\vb) = \frac{1}{(2\pi)^{d/2}} \e^{\left(-\frac{|\vb|^2}{2}\right)}\left(1+ \alpha \sum_{\ell=1}^d \cos(k_{\ell}x_{\ell})\right).
\end{equation}
In our experiments, we choose $k_{\ell}=0.5$ and for the weak Landau damping experiments we set $\alpha=0.01$. One can linearize the electric field around the Maxwellian equilibrium and get the linear solution for the electric field which is a good approximation if the parameter $\alpha$ is small. For the chosen parameter $k_{\ell}=0.5$ the damping rate of the electric field is $-0.1533$ according to the linear theory \cite{Sonnendruecker10} , i.e., the electric energy is damped by a factor $-0.3066$. We solve the weak Landau damping problem on the domain $[0,4\cdot \pi]^d \times [-6,6]^d$ discretized with a grid of $32$ points along each spacial dimension and $128$ points along the velocity dimensions. The experiment is done in one, two and three dimensions in TT format and in one and two dimensions on the full grid with the same resolution. The TT rounding is done to the accuracy $\varepsilon = 4 \cdot 10^{-6}$. Figure \ref{fig:linlan_eenergy} shows the electric energy as a function of time together with the envelope functions predicted by linear theory. We note that we recover the damping rate in all cases. We also see that the solution obtained in TT format is in good agreement with the solution on the full grid. Especially, the spurious recurrence occurs around time 63 on both the full and the TT grid. The maximum rank combinations are given in the upper part of Table \ref{tab:linlanCOMPRESSION} together with the corresponding compression rate compared to the full-grid solution.  Note that the ranks are only checked after each time step.

At time $t=40$, the $\ell_{\infty}$ error in the distribution function on the $32 \times 128$ grid in 2D is about $4.7 \cdot 10^{-6}$ compared to the solution on a refined grid.. In the tensor train computations with $\varepsilon=4 \cdot 10^{-6}$ for time $t=80$, we have an error of $2.0 \cdot 10^{-7}$ in the distribution function at time $t=40$. Hence, the TT truncation error is much smaller than the numerical error on the full grid for the chosen parameters.

\begin{table}
\footnotesize
\caption{Weak Landau damping. Compression on TT grid. Grid size: $32^d\times 128^d$.}
\begin{center}
\begin{tabular}{|c|c|c|c|}
\hline
initial value & & maximal rank combination & compression rate \\
\hline
\multirow{3}{*}{\eqref{eq:weaklan_aligned}}&2D & 7 & $2.7 \cdot 10^{-1}$ \\
&4D & 10,4,9 & $2.9 \cdot 10^{-4}$ \\
&6D & 11,4,18,4,10 & $2.5 \cdot 10^{-7}$ \\
\hline
\multirow{2}{*}{\eqref{eq:weaklan_diag}}&4D & 25,34,25 & $3.6 \cdot 10^{-3}$ \\
&6D & 21,23,23,23,18 & $1.7 \cdot 10^{-6}$ \\
\hline
\end{tabular}
\end{center}
\label{tab:linlanCOMPRESSION}
\end{table}

Admittedly, the problem is particularly suited for the tensor train format since the initial perturbation is aligned with the coordinate axes. We have therefore repeated the experiment with the initial condition 
\begin{equation}\label{eq:weaklan_diag}
f_0(\xb,\vb) = \frac{1}{(2\pi)^{d/2}} \e^{-\frac{|\vb|^2}{2}}\left(1+ \alpha  \cos\left(k\sum_{\ell=1}^dx_{\ell}\right)\right).
\end{equation}
For such a perturbation diagonal to the coordinate axes, the ranks between different pairs of $(x_i,v_i)$ are no longer small compared to the ranks between $x_i$ and $v_i$. The compression is reduced by approximately one order of magnitude as can be seen from the lower part of Table \ref{tab:linlanCOMPRESSION}. However, we still have a very good compression.

\subsection{Strong Landau damping}

If we increase the value of $\alpha$, linear theory is no longer a good approximation of the actual situation and nonlinear effects start to dominate. The filaments in the distribution function cannot be properly resolved on a given grid after a certain time. In 2D, the $\ell_{\infty}$ difference between the solution at time $t=30$ on a grid with $32\times 128$ points and a grid with $64\times 256$ is $1.2 \cdot 10^{-1}$. However, one is often only interested in some functional of the distribution function which can be recovered more accurately than the distribution function itself. In our example, we consider the electric field. The (absolute) $\ell_{\infty}$ error in the electric field is $3.0 \cdot 10^{-3}$ which corresponds to a relative error of about 10 \%. For one and two dimensions, we compare the electric field for the full and TT grid solution in Table \ref{tab:nonlinlanERROR}. For the given choice of the rounding threshold the error due to the rounding is comparable to the numerical error on the grid. Table \ref{tab:nonlinlanTIMES} reports the computing time and memory consumption for the representation of the distribution function. The reported times should give an indication on the order of magnitude of the computational time. In order to get quantitative results, a high-performance implementation of the methods needs to be used. Comparing the results with and without projection to conserve mass and momentum, we observe that there is no significant impact on the accuracy in the electric field nor on memory consumption or computing time. 

\begin{figure}
   \centering
      \subfloat[Electric energy.]{\includegraphics[scale=0.5]{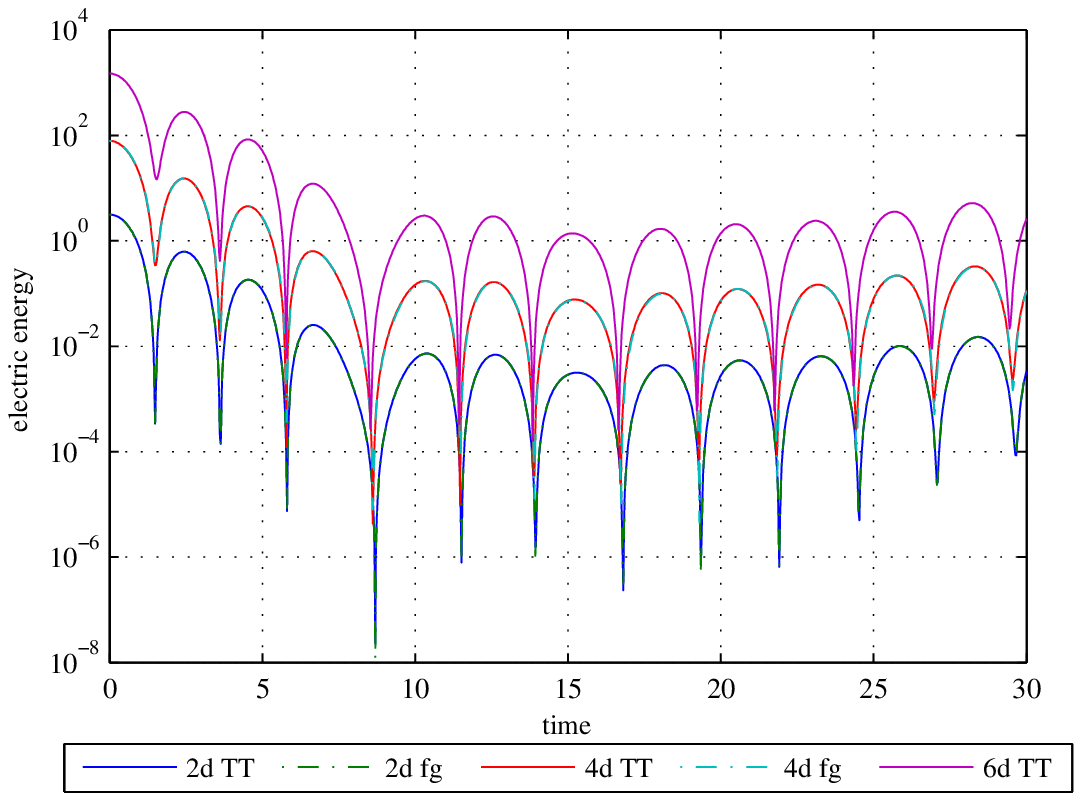}\label{fig:nonlinlanEEN}}\qquad
      \subfloat[Mass.]{\includegraphics[scale=0.5]{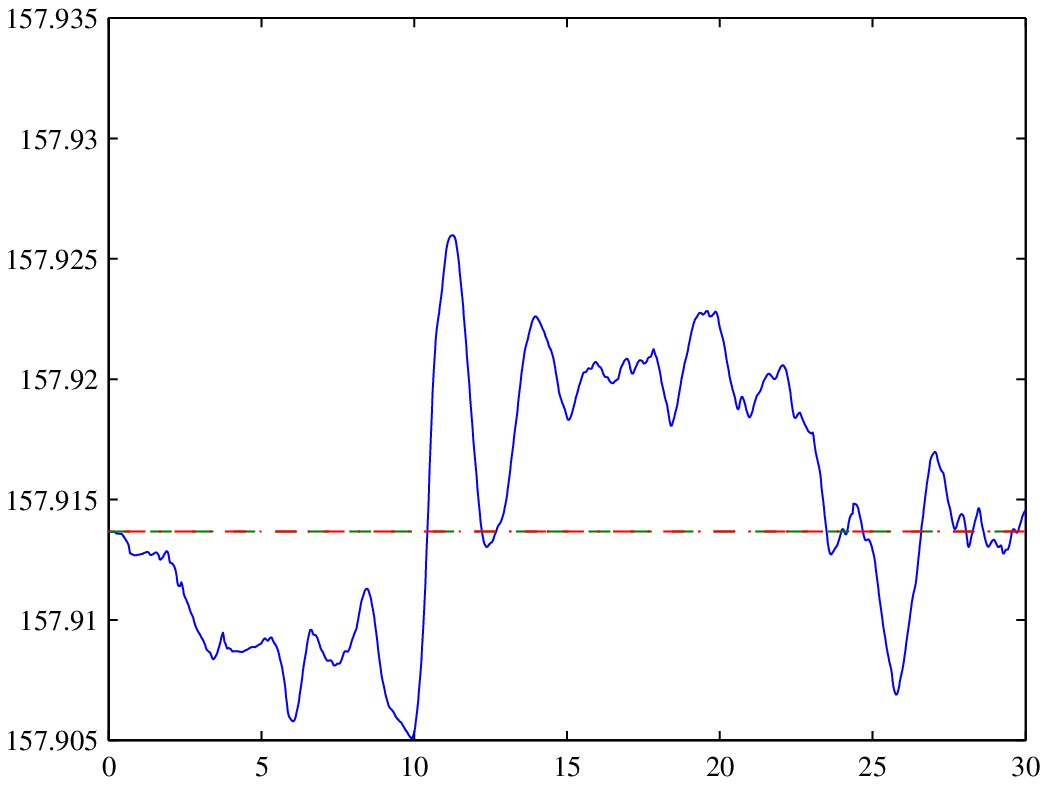}}\qquad
      \subfloat[Momentum, component 1.]{\includegraphics[scale=0.5]{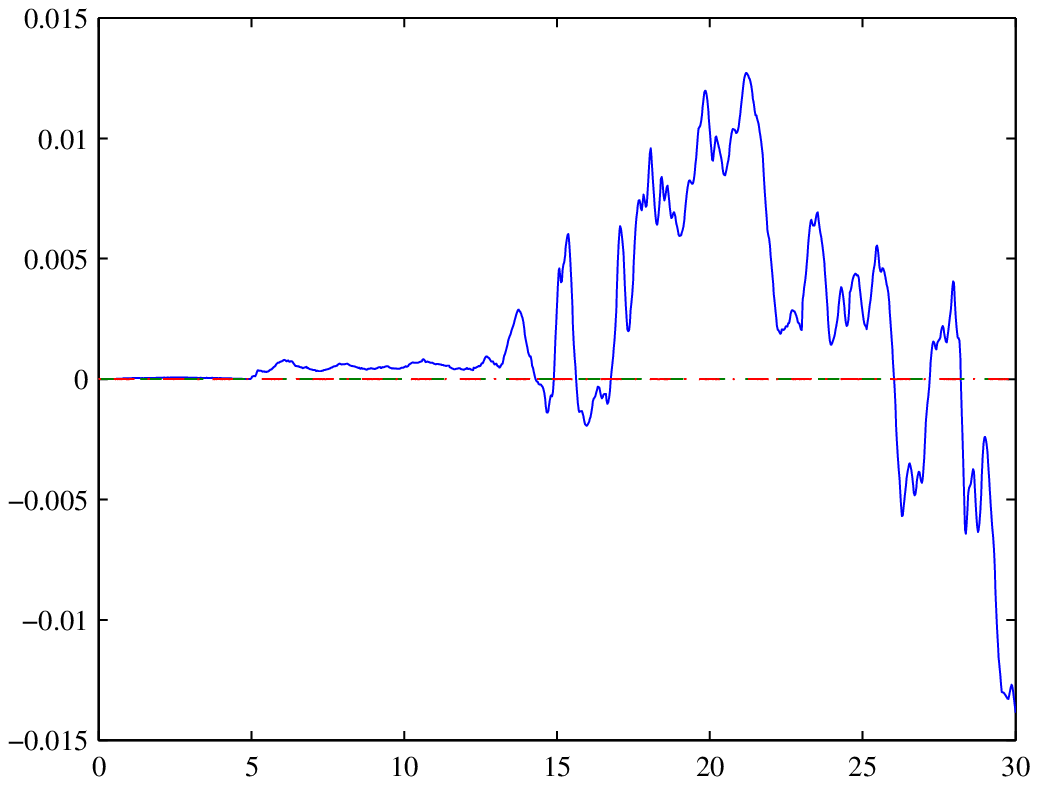}}\qquad
      \subfloat[Momentum, component 2.]{\includegraphics[scale=0.5]{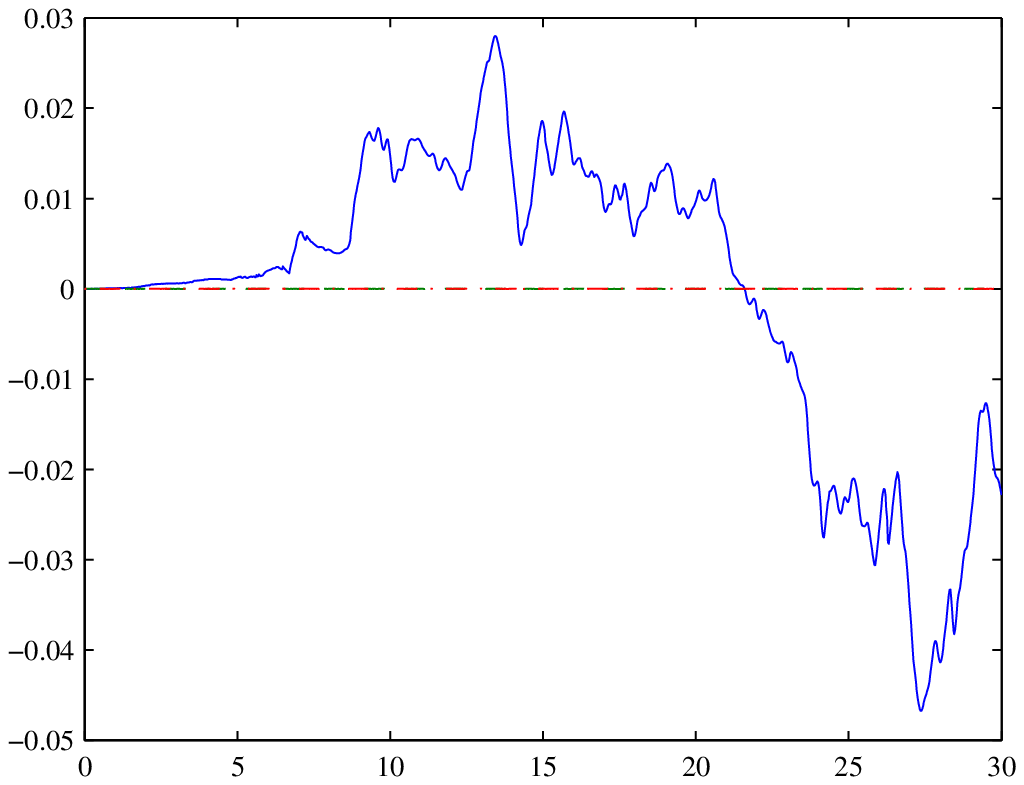}}\qquad
      \subfloat[$\ell_2$ norm.]{\includegraphics[scale=0.5]{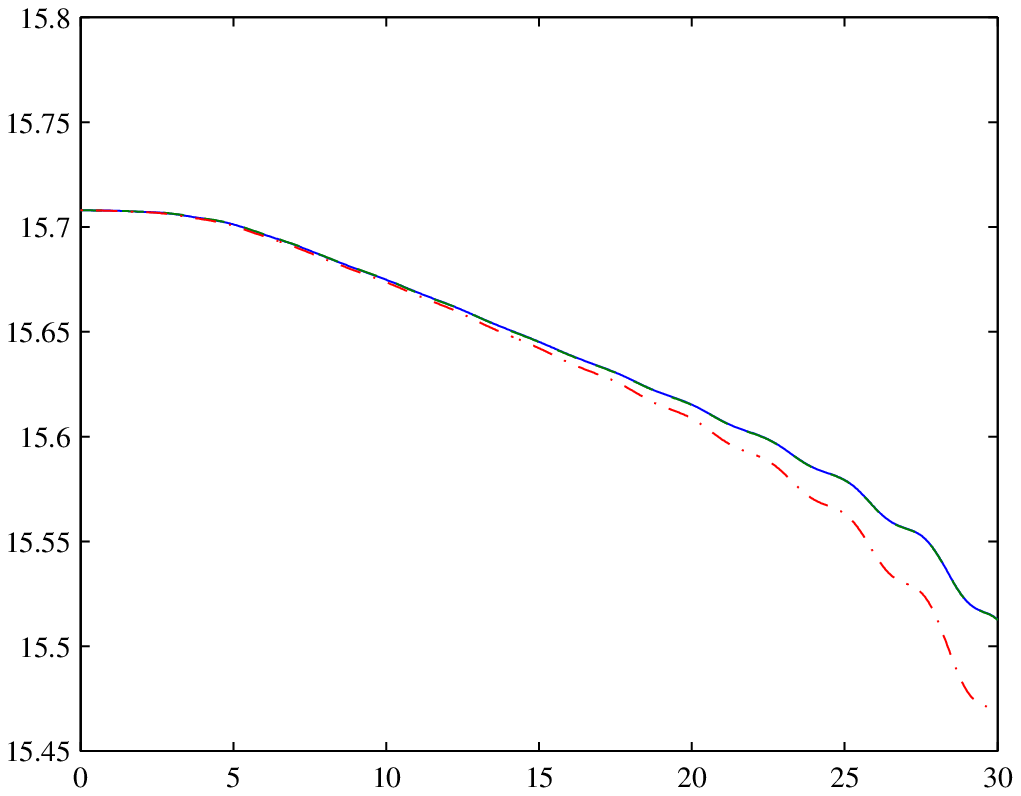}}\qquad
      \subfloat[Energy.]{\includegraphics[scale=0.5]{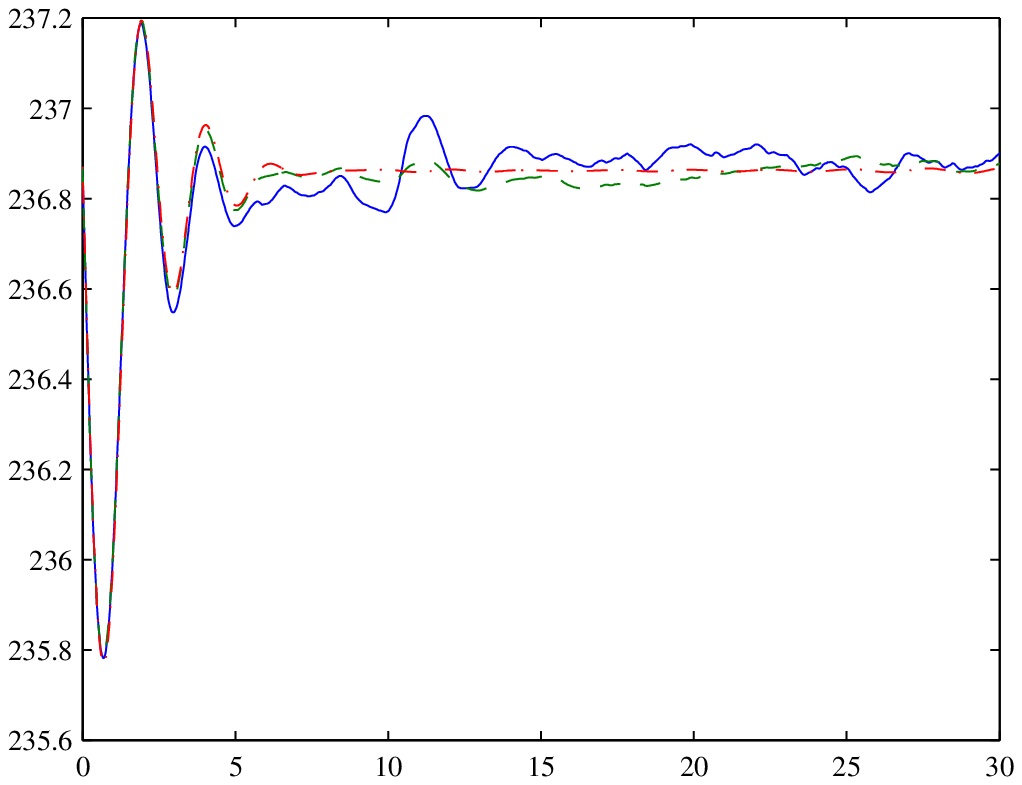}}\qquad
   \caption[Titel des Bildes]{Nonlinear Landau damping. (a) Electric energy for 2D, 4D and 6D computations. (b-f) Conservation of various properties for 4D simulations with ($- -$) and without ($-$) projection of mass and momentum compared to the full-grid simulation ($-\cdot$).}\label{fig:nonlinlan}
\end{figure}

In Figure \ref{fig:nonlinlanEEN}, the evolution of the electric energy is shown for the various runs. For the two dimensional problem, Figure \ref{fig:nonlinlan} shows the development over time of mass, momentum, $\ell_2$ norm, and energy and compares the version with and without projection to conserve mass and momentum. It can be seen that there is a considerable drift in mass and momentum if we do not project the solution. On the other hand, the figures show that the projection of mass and momentum nicely conserves these quantities without imparing the conservation of energy or $\ell_2$ norm.

\begin{table}[h]
\footnotesize
\caption{Strong Landau damping. Electric energy for the TT representation with (TTP) and without (TT) projection to conserve mass and momentum compared to the full grid. Grid size: $32^d\times 128^d$.}
\begin{center}
\begin{tabular}{|c|c|c|c|}
\hline
&$\varepsilon$& $\ell_{\infty}$ error $E_x$, ($E_y$)  \\
\hline\hline
2D TT & $4 \cdot 10^{-3}$ & $7.0 \cdot 10^{-4}$\\
2D TTP & $4 \cdot 10^{-3}$ & $6.4 \cdot 10^{-4}$\\
\hline
4D TT & $4 \cdot 10^{-4}$ & $2.4 \cdot 10^{-3}$, $2.3 \cdot 10^{-3}$\\
4D TTP & $4 \cdot 10^{-4}$ & $2.1 \cdot 10^{-3}$, $2.5 \cdot 10^{-3}$\\
\hline
\end{tabular}
\end{center}
\label{tab:nonlinlanERROR}
\end{table}

\begin{table}
\footnotesize
\caption{Strong Landau damping. Computing time (wall time) in seconds and memory of a TT representation with (TTP) and without (TT) projection to conserve mass and momentum compared to the solution on the full grid. Grid size: $32^2\times 128^2$, threshold for TT rounding: $4 \cdot 10^{-3}$ (2D), $4 \cdot 10^{-4}$ (4D), $2 \cdot 10^{-4}$ (6D).}
\begin{center}
\begin{tabular}{|c|c|cc|cc|}
\hline
dim & method & \# doubles for $f$ & fraction & wall time [s] & fraction \\
\hline
\hline
2D & FG & 4096 &  & $1.4 \cdot 10^1$ & \\
2D & TT & 2720 & $0.66 $  & $1.8 \cdot 10^1$& 1.3\\
2D & TTP& 3040 & $0.74 $ &  $2.0 \cdot 10^1$& 1.4\\ 
\hline
4D & FG & $1.7 \cdot 10^7$ &  & $6.2 \cdot 10^4$ & \\ 
4D & TT  & $5.9 \cdot 10^{4} $  & $3.5 \cdot 10^{-3}$ & $2.7 \cdot 10^2$ & $4.4 \cdot 10^{-3}$\\
4D & TTP & $6.0 \cdot 10^{4} $  & $3.6 \cdot 10^{-3}$ & $2.8 \cdot 10^2$ & $4.5 \cdot 10^{-3}$\\
\hline
6D & TT &  $7.1 \cdot 10^5$ & $1.0 \cdot 10^{-5}$ & $6.6 \cdot 10^3$ & \\
\hline
\end{tabular}
\end{center}
\label{tab:nonlinlanTIMES}
\end{table}

\subsection{Two stream instablity}

\begin{figure}[h]
   \centering
      \subfloat[Electric energy.]{\includegraphics[scale=0.5]{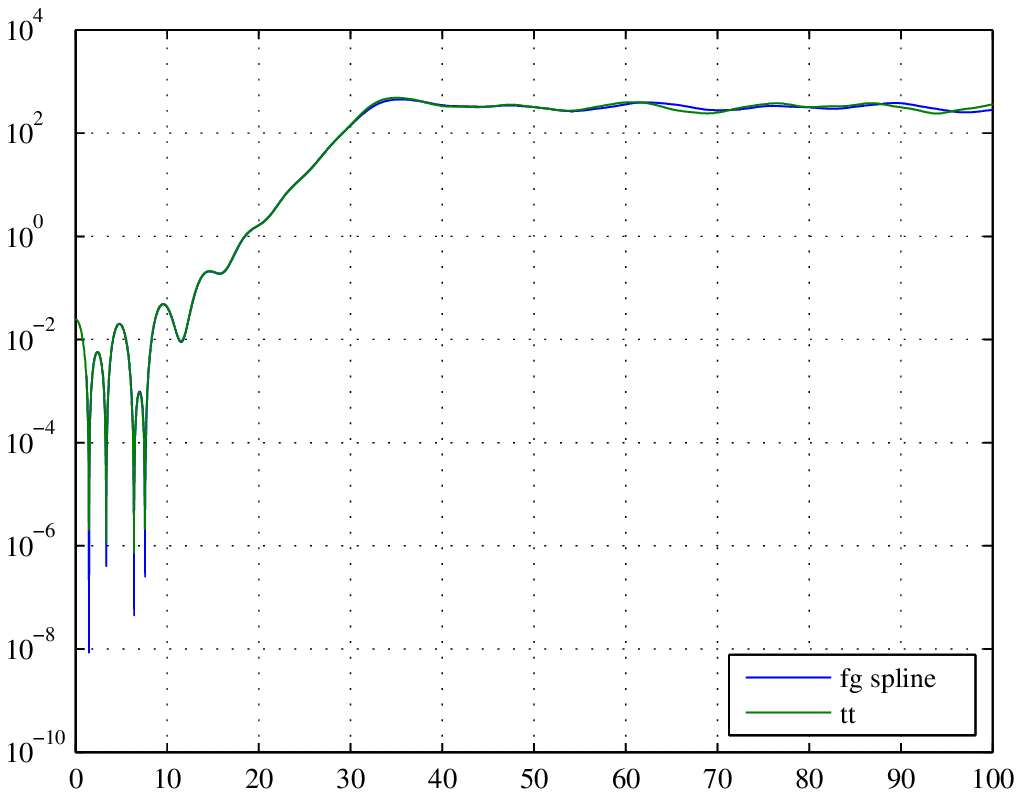}\label{fig:tsi_energy}}\qquad
      \subfloat[Ranks of TT representation.]{\includegraphics[scale=0.5]{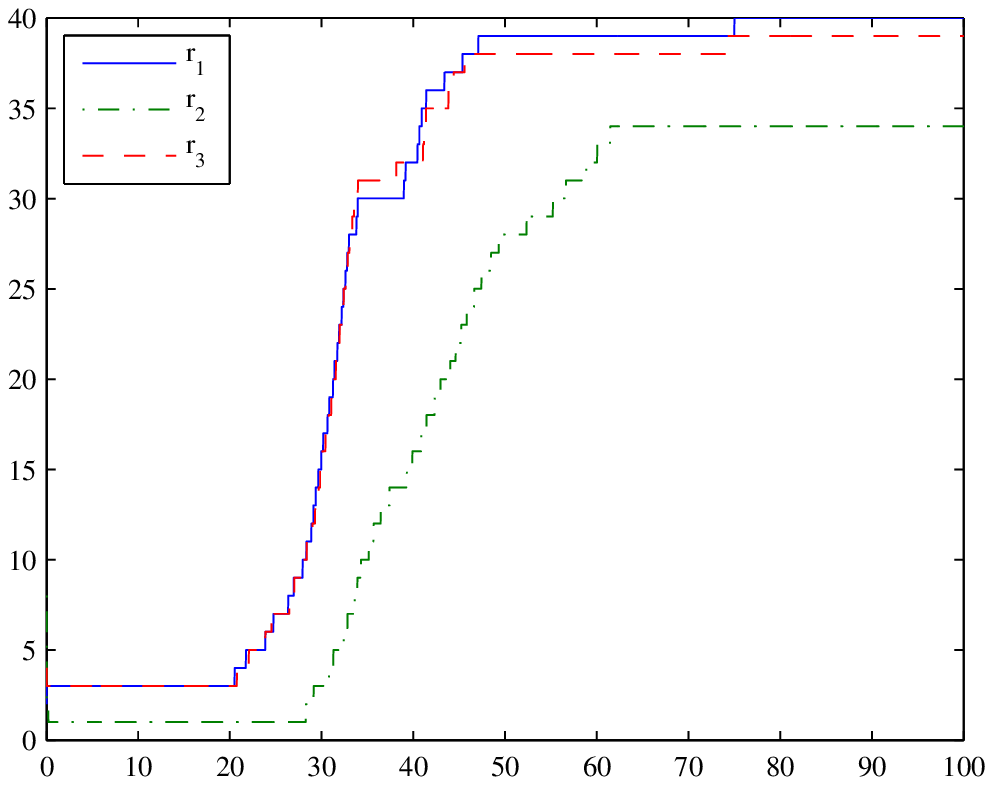}\label{fig:tsi_ranks}}\qquad
       \caption[Titel des Bildes]{Two stream instability.}\label{fig:twostream}
\end{figure}

Finally, we consider the two stream instability in 4D. In 2D phase space, the initial condition is
\begin{equation*}
f_0(x,v) = (1+\varepsilon\cos(k x))0.5 /\sqrt{2\pi}  \left(\e^{-0.5 (v-v_0)^2}+ \e^{-0.5(v+v_0)^2}\right).
\end{equation*}
In our simulations, we choose the parameters $k=0.2$, $\varepsilon=0.001$, and $v_0 = 2.4$.
We consider two kinds of extensions to 4D phase space
\begin{equation*} \begin{aligned}
f_0(\xb,\vb) =& \frac{0.5}{2 \pi} (1+\varepsilon\cos(k x_1))  \left(\e^{-0.5 (v_1-v_0)^2}+ \e^{-0.5(v_1+v_0)^2}\right)\e^{-0.5(v_2)^2},\\
f_0(\xb,\vb) =&  \frac{0.25}{2 \pi}(1+\varepsilon(\cos(k x_1)+ \cos(k x_2))) \left(\e^{-0.5 (v_1-v_0)^2}+ \e^{-0.5(v_1+v_0)^2}\right) \cdot\\
&\left(\e^{-0.5 (v_2-v_0)^2}+ \e^{-0.5(v_2+v_0)^2}\right),
\end{aligned}\end{equation*}
In the first case, we have an equilibrium state in $y, v_y$ plane. A solution in TT format detects this simple form and the solution is represented as a TT tensor with only rank $r_1$ different from one. The second case is a tensor product of two one dimensional two stream instabilities. Figure \ref{fig:tsi_ranks} shows the (inner) ranks as a function of time for a simulation on a grid with $64^4$ points and a rounding threshold $\varepsilon = 5 \cdot 10^{-4}$.
We see that the compression is very good in the beginning until about time 20. Thereafter the instability grows rapidly until about time 30. During this phase the ranks $r_1,r_3$ coupling the pairs $(v_x,x)$ and $(v_y,y)$ strongly increase. When nonlinear effects start to dominate and the electric energy flattens out, also rank $r_2$ increases for some time. Finally, the ranks remain almost constant from about time 60.

\begin{figure}[h]
   \centering
            \subfloat[tensor train solution]{\includegraphics[scale=0.5]{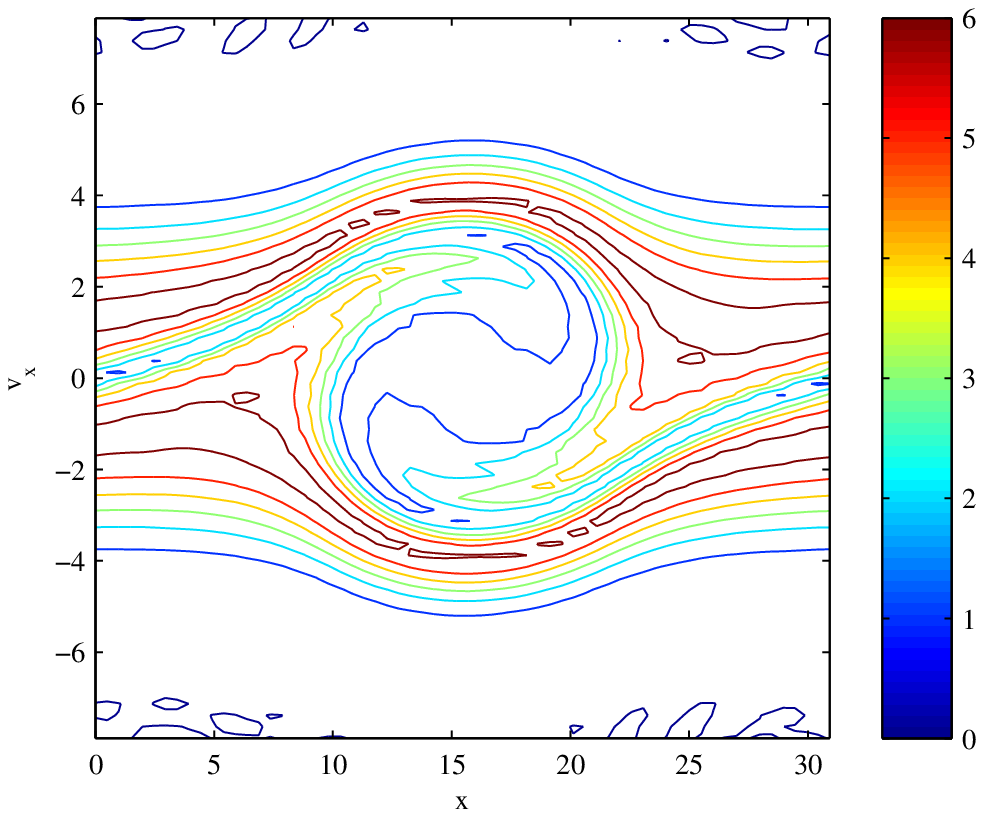}}\qquad
            \subfloat[full grid solution]{\includegraphics[scale=0.5]{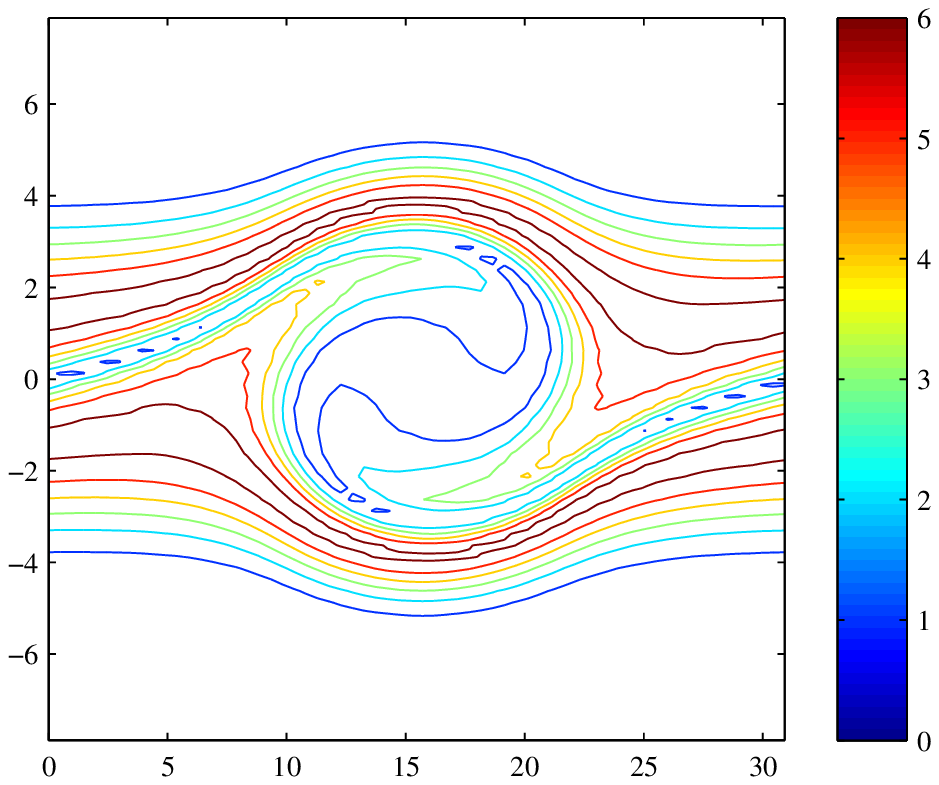}}\qquad
       \caption[Titel des Bildes]{$(v_x,x)$ projection of distribution function at time 36 for two stream instability.}\label{fig:twostream_phasespace}
\end{figure}

Figure \ref{fig:tsi_energy} shows the electric field as a function of time for the TT solution as well as the full grid solution. The curves show good agreement. In the nonlinear phase they start to deviate up to 27 \%. However, comparing the solution on the full grid with a solution on the same grid but with the same interpolation formulas as in the TT algorithm we see a deviation of up to 30 \%. Hence, the error due to TT rounding is on the scale of the numerical error. Figure \ref{fig:twostream_phasespace}  shows the distribution function in $(v_x,v)$ plane (integrated over $v_y,y$) at time 35 for the TT compressed as well as the full-grid solution with splines. We see that the TT solution covers the overall features of the full-grid solution but the solution is less smooth.

\section{Conclusions and Outlook}\label{sec:conclusions}

In this paper, we have devised a semi-La\-grangian Vlasov--Poisson solver with representation of the distribution function in tensor train format. For the efficient implementation of the advection step it is important to avoid direct matrix-vector products. Instead, we propose to compress the matrix describing the interpolation step and an iterative implementation of the arising Hadamard products. The method has been tested for a number of standard test cases in two to six dimensions. We have demonstrated that the solution can be compressed to a very high extent without losing essential parts of the solution when using a tensor train representation. Also the computing time is considerably reduced. As expected the gains from the tensor train representation become larger with growing dimension. 

In order to be able to study more complicated equations with more pronounced multidimensional effects, we plan to work on a high-performance implementation of the method. The choice of the interpolation formula and the effects of the CFL-like conditions and possible alleviations need to be studied in future work. Moreover, the effect of the rounding parameter and possibilities of automatic tolerance detection need further attention.

\section*{Acknowledgements}
The author thanks Eric Sonnendr{\"u}cker for bringing the tensor train framework to her attention and discussions on various aspects of this work. Discussions with Michel Mehrenberger and Marco Restelli on test cases were also appreciated.


\begin{thebibliography}{10}

\bibitem{Arber02}
T.~D. Arber and R.~G. Vann.
\newblock A critical comparison of eulerian-grid-based {V}lasov solvers.
\newblock {\em J. Comput. Phys.}, 180:339--357, 2002.

\bibitem{Besse04}
N.~Besse.
\newblock Convergence of a semi-{L}agrangian scheme for the one-dimensional
  {V}lasov--{P}oisson system.
\newblock {\em SIAM J. Numer. Anal.}, 42(1):350--382, 2004.

\bibitem{Buhmann00}
M.~D. Buhmann.
\newblock Radial basis functions.
\newblock {\em Acta Numerica}, 9:1--38, 2000.

\bibitem{Bungartz04}
H.-J. Bungartz and M.~Griebel.
\newblock Sparse grids.
\newblock {\em Acta Numerica}, 13:147--269, 2004.

\bibitem{Birdsall91}
A.~B.~L. C.~K.~Birdsall.
\newblock {\em Plasma Physics via Computer Simulations}.
\newblock Adam Hilger, 1991.

\bibitem{Cheng76}
C.~Z. Cheng and G.~Knorr.
\newblock The integration of the {V}lasov equation in configuration space.
\newblock {\em J. Comput. Phys.}, 22(3):330--351, 1976.

\bibitem{Cheng13}
Y.~Cheng, I.~M. Gamba, and P.~J. Morrison.
\newblock Study of conservation and recurrence of {R}unge--{K}utta
  discontinuous {G}alerkin schemes for {V}lasov--{P}oisson systems.
\newblock {\em J. Sci. Comput.}, 56(2):319--349, 2013.

\bibitem{Crouseilles08a}
N.~Crouseilles, M.~Gutnic, G.~Latu, and E.~Sonnendr{\"u}cker.
\newblock Comparison of two {E}ulerian solvers for the four-dimensional
  {V}lasov equation: Part {I}.
\newblock {\em Commun. Nonlinear Sci.}, 13(1):88--93, 2008.

\bibitem{Crouseilles08b}
N.~Crouseilles, M.~Gutnic, G.~Latu, and E.~Sonnendr{\"u}cker.
\newblock Comparison of two {E}ulerian solvers for the four-dimensional
  {V}lasov equation: Part {II}.
\newblock {\em Commun. Nonlinear Sci.}, 13(1):94--99, 2008.

\bibitem{Crouseilles10}
N.~Crouseilles, M.~Mehrenberger, and E.~Sonnendr{\"u}cker.
\newblock Conservative semi-{L}agrangian schemes for {V}lasov equations.
\newblock {\em J. Comput. Phys.}, 229(6):1927--1953, 2010.

\bibitem{Ayuso12}
B.~A. de~Dios, J.~A. Carrillo, and C.-W. Shu.
\newblock Discontinuous {G}alerkin methods for the multi-dimensional
  {V}lasov--{P}oisson problem.
\newblock {\em Math. Method Appl. Sci.}, 22(12):1250042, 2012.

\bibitem{Dolgov12a}
S.~Dolgov, B.~Khoromskij, and D.~Savostyanov.
\newblock Superfast {F}ourier transform using {QTT} approximation.
\newblock {\em Journal of Fourier Analysis and Applications}, 18(5):915--953,
  2012.

\bibitem{Dolgov12}
S.~V. Dolgov, B.~N. Khoromskij, and D.~V. Savostyanov.
\newblock Fast solution of parabolic problems in the tensor train/quantized
  tensor train format with initial application to the fokker-planck equation.
\newblock {\em SIAM J. Sci. Comput.}, 34:A2718--A2739, 2012.

\bibitem{Dolgov13}
S.~V. Dolgov and D.~V. Savostyanov.
\newblock Alternating minimal energy methods for linear systems in higher
  dimensions. {P}art {I}: {SPD} systems, 2013.

\bibitem{Dolgov13a}
S.~V. Dolgov and D.~V. Savostyanov.
\newblock Alternating minimal energy methods for linear systems in higher
  dimensions. {P}art {II}: {F}aster algorithm and application to nonsymmetric
  systems, 2013.

\bibitem{Dolgov14}
S.~V. Dolgov, A.~P. Smirnov, and E.~E. Tyrtyshnikov.
\newblock Low-rank approximation in the numerical modeling of the
  {F}arley--{B}uneman instability in ionospheric plasma.
\newblock {\em J. Comput. Phys.}, (0):--, 2014.

\bibitem{Filbet01}
F.~Filbet, E.~Sonnendr{\"u}cker, and P.~Bertrand.
\newblock Conservative numerical schemes for the {V}lasov equation.
\newblock {\em J. Comput. Phys.}, 172(1):166--187, 2001.

\bibitem{Grasedyck13}
L.~Grasedyck, D.~Kressner, and C.~Tobler.
\newblock A literature survey of low-rank tensor approximation techniques.
\newblock {arXiv} preprint 1302.7121v1, 2013.

\bibitem{Hackbusch12}
W.~Hackbusch.
\newblock {\em Tensor Spaces and Numerical Tensor Calculus}.
\newblock Springer Verlag, Berlin Heidelberg, 2012.

\bibitem{Hackbusch09}
W.~Hackbusch and S.~K{\"u}hn.
\newblock A new scheme for the tensor representation.
\newblock {\em J. Fourier Anal. Appl.}, 15:706--722, 2009.

\bibitem{Hairer06}
E.~Hairer, C.~Lubich, and G.~Wanner.
\newblock {\em Geometric numerical integration}.
\newblock Springer Verlag, Berlin Heidelberg, 2006.

\bibitem{Hatch12}
D.~Hatch, D.~del Castillo-Negrete, and P.~Terry.
\newblock Analysis and compression of six-dimensional gyrokinetic datasets
  using higher order singular value decomposition.
\newblock {\em J. Comput. Phys.}, 231(11):4234 -- 4256, 2012.

\bibitem{Holtz12}
S.~Holtz, T.~Rohwedder, and R.~Schneider.
\newblock The alternating linear scheme for tensor optimization in the tensor
  train format.
\newblock {\em SIAM J. Sci. Comput.}, 34:A683--A713, 2012.

\bibitem{Kazeev13}
V.~Kazeev, M.~Khammash, M.~Nip, and C.~Schwab.
\newblock Direct solution of the chemical master equation using quntized tensor
  trains.
\newblock Tech. Rep. 2013-04, Seminar for Applied Mathematics, ETH Z{\"u}rich,
  2013.

\bibitem{Kazeev12}
V.~Kazeev, O.~Reichmann, and C.~Schwab.
\newblock $hp$-{DG}-{QTT} solution of high-dimensional degenerate diffusion
  equations.
\newblock Tech. Rep. 2012-11, Seminar for Applied Mathematics, ETH Z{\"u}rich,
  2012.

\bibitem{Khoromskij12}
B.~N. Khoromskij.
\newblock Tensors-structured numerical methods in scientific computing: Survey
  on recent advances.
\newblock {\em Chemometr. Intell. Lab.}, 110(1):1--19, 2012.

\bibitem{Krall73}
N.~A. Krall and A.~W. Trivelpiece.
\newblock {\em Principles of Plasma Physics}.
\newblock McGrawHill, New York, 1973.

\bibitem{Lathauwer00}
L.~D. Lathauwer, B.~D. Moor, and J.~Vandewalle.
\newblock On the best rank-1 and rank-(r1 , ..., rn ) approximation of
  higher-order tensors.
\newblock {\em SIAM J. Matrix Anal. Appl.}, 21:1324--1342, 2000.

\bibitem{Mehrenberger13}
M.~Mehrenberger, C.~Steiner, L.~Marradi, N.~Crouseilles, E.~Sonnendr\"ucker,
  and B.~Afeyan.
\newblock {V}lasov on {GPU} ({VOG} project).
\newblock {\em ESAIM: Proc.}, 43:37--58, 2013.

\bibitem{Meyer90}
H.-D. Meyer, U.~Manthe, and L.~Cederbaum.
\newblock The multi-configurational time-dependent hartree approach.
\newblock {\em Chem. Phys. Lett.}, 165:73--78, 1990.

\bibitem{Oseledets11a}
I.~Oseledets.
\newblock {DMRG} approach to fast linear algebra in the {TT}-format.
\newblock {\em Comput. Methods Appl. Math.}, 11(3):272--403, 2011.

\bibitem{Oseledets11}
I.~Oseledets.
\newblock Tensor-train decomposition.
\newblock {\em SIAM J. Sci. Comput.}, 33(5):2295--2317, 2011.

\bibitem{Oseledets12}
I.~V. Oseledets, B.~N. Khoromskij, and R.~Schneider.
\newblock Efficient time-stepping scheme for dynamics on {T}{T}-manifolds.
\newblock Preprint~24, MPI MIS, 2012.

\bibitem{Qiu11}
J.~M. Qiu and C.~W. Shu.
\newblock Conservative semi-{L}agrangian finite difference {WENO} formulations
  with applications to the {V}lasov equation.
\newblock {\em Comm. Comput. Phys.}, 10, 2011.

\bibitem{Sonnendruecker10}
E.~Sonnendr\"ucker.
\newblock Approximation num{\'e}rique des {\'e}quations de {V}lasov--{M}axwell.
\newblock 2010.

\bibitem{Sonnendruecker99}
E.~Sonnendr{\"u}cker, J.~Roche, P.~Bertrand, and A.~Ghizzo.
\newblock The semi-{L}agrangian method for the numerical resolution of the
  {V}lasov equation.
\newblock {\em J. Comput. Phys.}, 149(2):201--220, 1999.

\bibitem{White92}
S.~R. White.
\newblock Density matrix formulation for quantum renormalization groups.
\newblock {\em Phys. Rev. Letters}, 69, 1992.

\end{thebibliography}
\end{document}